\documentclass[12pt]{article}
\usepackage{amsfonts}
\usepackage{hyperref}
\usepackage{mathrsfs}
\usepackage{indentfirst}
\usepackage{amsmath}
\usepackage{amssymb}
\usepackage[amsmath,thmmarks]{ntheorem}
\usepackage{graphicx}
\usepackage{tikz}

\newtheorem{definition}{\bf Definition}[section]
\newtheorem{lemma}{\bf Lemma}[section]
\newtheorem{theorem}{\bf Theorem}[section]
\newtheorem{remark}{\bf Remark}[section]

\newtheorem{example}{\bf Example}[section]

\newtheorem{proposition}{\bf Proposition}[section]
\newtheorem{comments}{\bf Comments}[section]

\begin{document}
\setcounter{page}{1}

\title{{\textbf{Characterizations of quasi-homogeneous aggregation functions}}\thanks {Supported by
the National Natural Science Foundation of China (No. 12071325)}}
\author{Feng-qing Zhu$^{1}$ \footnote{\emph{E-mail address}: zfqzfq2020@126.com} , Xue-ping Wang$^2$ \footnote{Corresponding author. xpwang1@hotmail.com}\\
\emph{\small {1.School of Sciences, Southwest Petroleum University, Chengdu 610500,}}\\\emph{\small {Sichuan, People's Republic of China}}\\
\emph{\small {2.School of Mathematical Sciences, Sichuan Normal University, Chengdu 610066,}}\\
\emph{\small {Sichuan, People's Republic of China}}}

\newcommand{\pp}[2]{\frac{\partial #1}{\partial #2}}
\date{}
\maketitle
\begin{quote}
{\bf Abstract}  In this article, we first give the characterizations of quasi-homogeneous aggregation functions, which show us that quasi-homogeneous aggregation functions are classified into three classes. We then introduce the concept of triple generator of quasi-homogeneous aggregation function, which is applied to construct a quasi-homogeneous aggregation function.

\end{quote}

\section{Introduction}

The homogeneous functions play an important role in image processing, decision making and other relevant fields, see, e.g., \cite{Cabrera2014,Jurio2013, Jurio2014}. Thus, it is very valuable to study them from the theoretical point of view. A function $A:[0, 1]^2\rightarrow[0, 1]$ is said to be homogeneous of order $k>0$ if it satisfies $A(\lambda x, \lambda y) =\lambda^{k}A(x, y)$ for all $ x, y, \lambda\in[0, 1]$. So far, there are many articles for investigating the homogeneity of particular aggregation functions. For example, R\"{u}ckschlossov\'{a} presented a complete characterization of homogeneous aggregation functions \cite{T2005}. On the other hand, a more relaxed homogeneity, called a quasi-homogeneity, was introduced by Ebanks in \cite{Ebanks1998}. Quasi-homogeneous t-norms are defined by $T(\lambda x, \lambda y) =\varphi^{-1}(\psi(\lambda)\varphi(T(x, y)))$ for all $ x, y, \lambda\in[0, 1]$ where $\psi:[0, 1] \rightarrow[0, 1]$ is a function and $\varphi:[0, 1] \rightarrow[0, \infty)$ is a continuous injection and $T$ is a triangular norm \cite {Ebanks1998}. Just replacing the triangular norm $T$ by a copula, quasi-homogeneous copulas were similarly defined by Mayor, Mesiar and Torrens \cite{Mayor2008}. Su, Zong and Mesiar \cite{Su2022} investigated the characterizations of homogeneous and quasi-homogeneous aggregation functions, respectively. Recently, Wang and Zhu \cite{Zhu} studied the pseudo-homogeneous overlap and grouping functions. One can easily see that all the results of \cite{Su2022,Zhu} are suitable for an aggregation function whose diagonal function is continuous. So that a natural problem is: what are the characterizations of quasi-homogeneous aggregation functions? This article will pay attention to this problem.

The rest of this article are organized as follows. In Section 2, we give some comments on Proposition 5 and Theorem 6 in \cite{Su2022}. In Section 3, we show the characterizations of quasi-homogeneous aggregation functions. In Section 4, we introduce the concept of triple generator of quasi-homogeneous aggregation function, which is applied to construct a quasi-homogeneous aggregation function. A conclusion is drawn in Section 5.

\section{Comments on some results presented in \cite{Su2022}}

In this section, we first modify the definition of quasi-homogeneous aggregation function introduced by \cite {Su2022}, and then give some comments on Proposition $5$ and Theorem $6$ obtained by \cite{Su2022}.

Firstly, recall that a function $A: [0, 1]^2\rightarrow[0, 1]$ is called an \emph{aggregation function} if it is increasing and satisfies the boundary conditions $A(0, 0) = 0$ and $A(1, 1) = 1$.  A function $\delta_{A}: [0, 1] \rightarrow[0, 1]$ with
$\delta_{A}(x) = A(x, x)$ is called a \emph{diagonal function} of $A$, in symbols $\delta_A$.

\begin{definition}\label{def2.2}(\cite{Su2022})
\emph{An aggregation function $A:[0, 1]^2\rightarrow[0, 1]$ is said to be $(\varphi, \psi)$-quasi-homogeneous
if, for all $ x, y, \lambda\in[0, 1]$,
 \begin{equation*}
 A(\lambda x, \lambda y) =\varphi^{-1}(\psi(\lambda)\varphi(A(x, y)))
  \end{equation*}
where $\psi:[0, 1] \rightarrow[0, 1]$ is an arbitrary function and $\varphi:[0, 1] \rightarrow[0, \infty)$ is a continuous injection.}
\end{definition}

It is well known that a continuous injection maps a compact set to a compact set. Therefore, to ensure the existence of the inverse function $\varphi^{-1}$ of function $\varphi$ since $[0, \infty)$ is not compact, from the mathematical point of view, we should reconsider the mapping $\varphi$. This leads us to slightly modify the definition of $(\varphi, \psi)$-quasi-homogeneous aggregation functions as follows.
\begin{definition}\label{def2.1}
\emph{An aggregation function $A:[0, 1]^2\rightarrow[0, 1]$ is said to be} $(\varphi, \psi)$-quasi-homogeneous
(quasi-homogeneous \emph{for short) if, for all $ x, y, \lambda\in[0, 1]$,}
 \begin{equation}\label{eq2.1}
 A(\lambda x, \lambda y) =\varphi^{-1}(\psi(\lambda)\varphi(A(x, y)))
  \end{equation}
\emph{where $\psi:[0, 1] \rightarrow[0, 1]$ is a function and $\varphi:[0, 1] \rightarrow[0, b]$ is a continuous bijection with $[0, b]\subseteq[0, \infty]$.}
\end{definition}

The following are Proposition $5$ and Theorem $6$ of \cite{Su2022}, respectively.
\begin{lemma}\label{Lemma2.1}(see Proposition 5 in \cite{Su2022})
If an aggregation function $O$ is $(\varphi, \psi)$-quasi-homogeneous, then
\begin{enumerate}
\item [(\romannumeral 1)] its diagonal section $\delta_O(x)=O(x, x)$ is an increasing bijection.
\item [(\romannumeral 2)] $\psi(x)=x^c$ and $\varphi(x)=(\delta^{-1}_O(x))^c$ for some arbitrarily chosen $c>0$.
\item [(\romannumeral 3)] $\delta^{-1}_O\circ O$ is homogeneous of order 1.
\end{enumerate}
\end{lemma}

\begin{lemma}\label{Lemma2.2}(see Theorem 6 in \cite{Su2022})
An aggregation function $O$ is $(\varphi, \psi)$-quasi-homogeneous if and only if its diagonal $\delta_O$ is an increasing bijection and there exist increasing functions $h, g:[0, 1] \rightarrow[0, 1]$ fulfilling that $h(1)=g(1)=1$, $\frac {h(x)}{x}$ and $\frac {g(x)}{x}$ are decreasing on $(0, 1]$ such that
\begin{equation*}
O(x,y)=\left\{\begin{array}{ll}
0 & {\mbox{\scriptsize\normalsize if  } (x,y)=(0,0)},\\
\delta_O(yh(\frac{x}{y})) & {\mbox{\scriptsize\normalsize if  } x\leq y \ and\  y\neq0},\\
\delta_O(xg(\frac{y}{x})) & {\mbox{\scriptsize\normalsize if  } y\leq x \ and\ x\neq0}.
\end{array}
\right.
\end{equation*}
In this case, $O$ is $(\varphi, \psi)$-quasi-homogeneous with $\psi(x)=x^c$ and $\varphi(x)=(\delta_O^{-1}(x))^c$ for some arbitrarily chosen $c>0$.
\end{lemma}

To give comments on Lemmas \ref{Lemma2.1} and \ref{Lemma2.2}, we recall the following important result.
\begin{lemma}[see Theorem 13.1.6 in \cite{Kuczma1985}]\label{lem2.1}
Let $D$ be one of the sets $(0, 1)$, $[0, 1)$, $(-1, 1)$, $(-1, 0)\cup (0, 1)$, $(1, \infty)$, $(0, \infty)$, $[0, \infty)$, $(-\infty, 0)\cup (0, \infty)$ and $R$. A function $f:D\rightarrow R$ is a continuous solution of the multiplicative Cauchy equation $f(xy) = f(x)f(y)$ if and only if either $f=0$, or $f=1$, or $f$ has one of the following forms:
$$f(x)= |x|^c, x\in D,$$
$$f (x) = |x|^c \emph{sgn}(x), x\in D,$$
with a certain $c\in R$. If $0\in D$, then $c>0$.
\end{lemma}

\begin{comments}\label{Remark 2.1}
\emph{$(\romannumeral 1)$ To avoid a symbol confusion, in Lemmas \ref{Lemma2.1} and \ref{Lemma2.2} we use $\psi$ instead of $f$ in Proposition 5 and Theorem 6 of \cite{Su2022}, respectively.}

\emph{$(\romannumeral 2)$ It is well known that the aggregation function $T_{D}$ is given by
\begin{equation*}
T_{D}(x,y)=\left\{\begin{array}{ll}
0 & {\mbox{\scriptsize\normalsize if  } x,y \in[0, 1)},\\
{\mbox{min}\{x,y\}} & {\mbox{\scriptsize\normalsize otherwise}.}
\end{array}
\right.
\end{equation*} From Definition \ref {def2.1}, one easily verify that $T_{D}$ is a $(\varphi, \psi)$-quasi-homogeneous aggregation function, where  $$\psi(x)=\left\{\begin{array}{ll}
0 & {\mbox{\scriptsize\normalsize if } x\in[0, 1),}\\
1 & {\mbox{\scriptsize\normalsize if } x=1}
\end{array}
\right.$$
and $\varphi:[0, 1] \rightarrow[0, 1]$ is an increasing bijection. However, it is evident that the diagonal of $T_{D}$ is not an increasing bijection.
This indicates that Lemma \ref{Lemma2.1} (Proposition 5 in \cite{Su2022}) is incorrect. Unfortunately, in \cite{Su2022}, Theorem 6 is directly derived from Proposition 5. Therefore, Theorem 6 of \cite{Su2022} is incorrect, or more precisely, Theorem 6 of \cite{Su2022} is false for accurately characterizing the $(\varphi, \psi)$-quasi-homogeneous aggregation function}

\emph{$(\romannumeral 3)$ It is clearly that Lemma \ref{lem2.1} (Theorem 6 in \cite{Su2022}) is not suitable for $D=[0,1]$. Unfortunately, the authors happened to directly obtain the expression of function $f$ by applying Theorem 13.1.6 in \cite{Kuczma1985} when they proved Proposition 5 in \cite{Su2022}.}
\end{comments}

Next, we provide a complete characterization of quasi-homogeneous aggregation functions.

\section{Characterizations of quasi-homogeneous aggregation functions}

This section is devoted to characterize quasi-homogeneous aggregation functions.
\begin{lemma}\label{lem3.1}
If an aggregation function $A:[0, 1]^2\rightarrow[0, 1]$ is $(\varphi, \psi)$-quasi-homogeneous, then
\begin{enumerate}
\item [(\romannumeral 1)] $\varphi(1)\neq 0$;
\item [(\romannumeral 2)] $\varphi(0)=0$;
\item [(\romannumeral 3)] $\varphi$ is an increasing bijection.
\end{enumerate}
\end{lemma}
\begin{proof}
$(\romannumeral 1)$ Suppose $\varphi(1)=0$. Then by Eq.(\ref{eq2.1}),
$$A(\lambda,\lambda) =\varphi^{-1}(\psi(\lambda)\varphi(A(1,1))) = \varphi^{-1}(0) = 1$$ for all $\lambda\in[0, 1],$
which is impossible. Therefore, $\varphi(1)\neq 0$.

$(\romannumeral 2)$ By Eq.(\ref{eq2.1}),
$$\delta_{A}(\lambda x)=A(\lambda x, \lambda x) =\varphi^{-1}(\psi(\lambda)\varphi(A(x, x)))=\varphi^{-1}(\psi(\lambda)\varphi(\delta_{A}(x)))$$
for all $ x, \lambda\in[0, 1]$, i.e., $\delta_{A}(\lambda x)=\varphi^{-1}(\psi(\lambda)\varphi(\delta_{A}(x)))$, or equivalently,
\begin{equation}\label{eq3.1}
\varphi(\delta_{A}(\lambda x))=\psi(\lambda)\varphi(\delta_{A}(x))
 \end{equation}
 for all $x, \lambda \in[0, 1]$. Considering $x=0$ in Eq.(\ref{eq3.1}), we have
$$\varphi(0)=\psi(\lambda)\varphi(0)\mbox{ for all } \lambda \in[0, 1],$$
or equivalently,
$$\varphi(0)(1-\psi(\lambda))=0\mbox{ for all } \lambda \in[0, 1],$$
then $$\varphi(0)=0.$$
Otherwise, $\psi(\lambda)=1\mbox{ for all }\lambda \in[0, 1]$. In this case, considering $\lambda=0$ and Eq.(\ref{eq2.1}), we have
$$A(0,0)=A(x,y)\mbox{ for all } x, y \in[0, 1],$$
a contradiction. Thus
$$\varphi(0)=0.$$

$(\romannumeral 3)$ Obviously, from $(\romannumeral 1)$ and $(\romannumeral 2)$ we have that $\varphi$ is an increasing bijection.
\end{proof}

By Lemma \ref{lem3.1}, if $A:[0, 1]^2\rightarrow[0, 1]$ is a $(\varphi, \psi)$-quasi-homogeneous aggregation function then $\varphi(1)=b$. In what follows, we shall distinguish $b$ by two cases $0<b<\infty$ and $b=\infty$, respectively.
\subsection{The case $0<b<\infty$}
\begin{lemma}[see Theorem 13.1.9 in \cite{Kuczma1985}]\label{lem2.2}
Let $D$ be one of the sets $(0, 1)$, $[0, 1)$, $(-1, 1)$, $(-1, 0)\cup (0, 1)$, $(1, \infty)$, $(0, \infty)$, $[0, \infty)$, $(-\infty, 0)\cup (0, \infty)$ and $R$, and let $f:D\rightarrow R$ is a solution of the multiplicative Cauchy equation $f(xy) = f(x)f(y)$. If $f$ is measurable, then it is continuous in $D\setminus\{0\}$.
\end{lemma}
\begin{lemma}\label{lem3.2}
If an aggregation function $A:[0, 1]^2\rightarrow[0, 1]$ is $(\varphi, \psi)$-quasi-homogeneous, then
\begin{enumerate}
\item [(\romannumeral 1)] $\psi$ is increasing;
\item [(\romannumeral 2)] $\psi(\lambda x)=\psi(\lambda)\psi(x)$ for all $\lambda, x\in[0, 1]$;
\item [(\romannumeral 3)] $\psi(0)=0$, $\psi(1)=1$;
\item [(\romannumeral 4)] $\psi$ is continuous on $(0,1)$.
\end{enumerate}
\end{lemma}
\begin{proof}
$(\romannumeral 1)$ Considering $x=1$ in Eq.(\ref{eq3.1}), then
\begin{equation}\label{eq3.3}
\psi(\lambda)=\frac{1}{\varphi(1)}\varphi(\delta_{A}(\lambda))\mbox{ for all } \lambda \in[0, 1].
 \end{equation}
Lemma \ref{lem3.1} and Eq.(\ref{eq3.3}) imply that $\psi$ is increasing.\\
$(\romannumeral 2)$ By Eqs.(\ref{eq3.1}) and (\ref{eq3.3}), we have
 \begin{equation}\label{eq3.4}
\psi(\lambda x)=\psi(\lambda)\psi(x)\mbox{ for all } \lambda, x\in[0, 1].
 \end{equation}
$(\romannumeral 3)$ By Eq.(\ref{eq3.4}), we have
\begin{equation}\label{eq3.5}
\psi(0)=\psi(0)\psi(0)
\end{equation}
and
\begin{equation}\label{eq3.6}
\psi(1)=\psi(1)\psi(1).
\end{equation}
Eq.(\ref {eq3.5}) implies $\psi(0)=0$ or $\psi(0)=1$,
and
Eq.(\ref{eq3.6}) means $\psi(1)=0$ or $\psi(1)=1$. Then, $0\leq\psi(\lambda)\leq1$ for all $\lambda\in[0, 1]$ since $\psi$ is increasing.
If $\psi(0)=1$, then $\psi(\lambda)\equiv1$ for all $\lambda\in[0, 1]$ since $\psi$ is increasing. In this case, $A(0, 0)=\varphi^{-1}(\psi(0)\varphi(A(1, 1)))=1$, a contradiction.
Therefore, $\psi(0)=0$. Similarly, we can get $\psi(1)=1$.\\
$(\romannumeral 4)$
Using Lemma \ref{lem2.2}, it follows immediately from $(\romannumeral 1)$ and $(\romannumeral 2)$.
\end{proof}

\begin{lemma}\label{lemma3.1}
If an aggregation function $A:[0, 1]^2\rightarrow[0, 1]$ is $(\varphi, \psi)$-quasi-homogeneous, then one of the following statements holds:
\begin{enumerate}
\item [(\romannumeral 1)] $\psi(x)=x^c $ for a certain $c >0$.
\item [(\romannumeral 2)]
\begin{equation*}
\psi(x)=\left\{\begin{array}{ll}
0 & {\mbox{\scriptsize\normalsize if  } x \in[0, 1),}\\
1 & {\mbox{\scriptsize\normalsize if  } x=1.}
\end{array}
\right.
\end{equation*}
\item [(\romannumeral 3)] \begin{equation*}
\psi(x)=\left\{\begin{array}{ll}
0 & {\mbox{\scriptsize\normalsize if } x=0,}\\
1 & {\mbox{\scriptsize\normalsize if } x \in(0, 1].}
\end{array}
\right.
\end{equation*}
\end{enumerate}
\end{lemma}
\begin{proof}
By Lemmas \ref{lem2.1} and \ref{lem3.2}, we have that for any $x\in (0, 1)$, $\psi(x)=0$ or $\psi(x)=1$ or, $\psi(x)$ has one of the following forms:
$\psi(x)= |x|^c, $ $ \psi(x)= |x|^c \mbox{sgn}(x)$,
with a certain $c\in R$.
Next, we investigate the function $\psi$ by distinguishing three cases.

$(\romannumeral 1)$ If for any $x\in (0, 1)$, $\psi(x)$ has one of the following forms:
$\psi(x)= |x|^c, $ $ \psi(x)= |x|^c \mbox{sgn}(x)$,
with a certain $c\in R$, then by Lemma $\ref{lem3.2}$, $\psi(x)= x^c$ with a certain $c>0$.

$(\romannumeral 2)$ If $\psi(x)=0$ for any $x\in (0, 1)$, then, by $(\romannumeral 3)$ and $(\romannumeral 4)$ of Lemma \ref{lem3.2}, $$\psi(x)=\left\{\begin{array}{ll}
0 & {\mbox{\scriptsize\normalsize if } x \in[0, 1),}\\
1 & {\mbox{\scriptsize\normalsize if } x=1.}
\end{array}
\right. $$

$(\romannumeral 3)$ If $\psi(x)=1$ for any $x\in (0, 1)$, then, by $(\romannumeral 3)$ and $(\romannumeral 4)$ of Lemma \ref{lem3.2}, we have $$\psi(x)=\left\{\begin{array}{ll}
0 & {\mbox{\scriptsize\normalsize if } x=0,}\\
1 & {\mbox{\scriptsize\normalsize if } x \in(0, 1].}
\end{array}
\right. $$
\end{proof}

\begin{lemma}\label{lem3.3}
If an aggregation function $A:[0, 1]^2\rightarrow[0, 1]$ is $(\varphi, \psi)$-quasi-homogeneous with $\psi(x)=x^c$ for a certain $c>0$, then the following statements holds:
\begin{enumerate}
\item [(\romannumeral 1)] $\delta_{A}$ is an increasing bijection.
\item [(\romannumeral 2)] $\varphi(x)=\varphi(1)(\delta_{A}^{-1}(x))^c$.
\end{enumerate}
\end{lemma}
\begin{proof}
$(\romannumeral 1)$ From Eq.(\ref{eq2.1}) we have for any $x\in [0, 1]$,
$$A(x, x)=\varphi^{-1}(x^c\varphi(A(1,1)))=\varphi^{-1}(x^c\varphi(1))\mbox{ for a certain } c>0,$$
then for any $x\in [0, 1]$,
\begin{equation}\label{eq3.9}
\delta_{A}(x)=\varphi^{-1}(x^c\varphi(1))\mbox{ for a certain } c>0.
\end{equation}
Obviously, $\delta_{A}$ is an increasing bijection since $\varphi$ is an increasing bijection.

$(\romannumeral 2)$ Using Eq.(\ref{eq3.9}), we obtain that $\varphi(x)=\varphi(1)(\delta_{A}^{-1}(x))^c.$
\end{proof}

\begin{remark}\label{Remark3.1}
\emph{An aggregation function $A:[0, 1]^2\rightarrow[0, 1]$ is $(\varphi, \psi)$-quasi-homogeneous with $\psi(x)=x^c$ and $\varphi(x)=\varphi(1)(\delta_{A}^{-1}(x))^c$ for a certain $c>0$ if and only if it is $(\widetilde{\varphi}, \widetilde{\psi})$-quasi-homogeneous with $\widetilde{\psi}(x)=x^\alpha$ and $\widetilde{\varphi}(x)=(\delta_{A}^{-1}(x))^\alpha$ for any $\alpha>0$.}

    \emph{Indeed, for any }$x,y\in[0,1]$, $$\varphi^{-1}(\psi(\lambda)\varphi(A(x,y)))=\delta_{A}((\frac{\lambda^c \varphi(1)(\delta_{A}^{-1}(A(x,y)))^c}{\varphi(1)} )^{\frac{1}{c}})=\delta_{A}(\lambda \delta_{A}^{-1}(A(x,y)))$$ \emph{ and }
$$\widetilde{\varphi}^{-1}(\widetilde{\psi}(\lambda)\widetilde{\varphi}(A(x,y)))=\delta_{A}((\lambda^{\alpha} (\delta_{A}^{-1}(A(x,y)))^{\alpha} )^{\frac{1}{\alpha}})=\delta_{A}(\lambda \delta_{A}^{-1}(A(x,y))).$$
\end{remark}

Therefore, by Remark \ref{Remark3.1} and Lemma \ref{lem3.3}, in what follows, when we discuss the $(\varphi, \psi)$-quasi-homogeneity of aggregation function $A$ with $\psi(x)=x^c$ and $\varphi(x)=\varphi(1)(\delta_{A}^{-1}(x))^c$ for a certain $c>0$, we can always presuppose $\psi(x)=x$ and $\varphi(x)=\delta_{A}^{-1}(x)$ for any $x\in [0, 1]$.

\begin{theorem}\label{th3.2}
An aggregation function $A:[0, 1]^2\rightarrow[0, 1]$ is $(\varphi, \psi)$ \emph{-quasi-homogeneous} if and only if one of the following statements holds:
\begin{enumerate}
\item [(\romannumeral 1)] The diagonal $\delta_{A}$ is an increasing bijection and there exist two increasing functions $h, g : [0, 1] \rightarrow[0, 1]$ fulfilling that $h(1)=g(1)=1$, $\frac{\delta_{A}^{-1}(h(x))}{x}$ and $\frac{\delta_{A}^{-1}(g(x))}{x}$ are decreasing on $(0, 1]$ such that for all $x, y\in [0,1]$
\begin{equation}\label{eq3.10}
A(x,y)=\left\{\begin{array}{ll}
0 & {\mbox{\scriptsize\normalsize if  } (x,y)=(0,0)},\\
\delta_{A}(y\delta_{A}^{-1}(h(\frac{x}{y}))) & {\mbox{\scriptsize\normalsize if  } x\leq y \ and\  y\neq0},\\
\delta_{A}(x\delta_{A}^{-1}(g(\frac{y}{x}))) & {\mbox{\scriptsize\normalsize if  } y\leq x \ and\ x\neq0},
\end{array}
\right.
\end{equation}
$\psi(x)=x$ and $\varphi(x)=\delta_{A}^{-1}(x)$.
\item [(\romannumeral 2)] There exist two constants $\alpha, \beta \in [0, 1]$ such that for all $x, y\in [0,1]$
\begin{equation}\label{eq3.11}
A(x,y)=\left\{\begin{array}{ll}
0 & {\mbox{\scriptsize\normalsize if  } (x,y)=(0,0)},\\
1 & {\mbox{\scriptsize\normalsize if  } x,y \in(0, 1]},\\
\alpha & {\mbox{\scriptsize\normalsize if  } x=0, y\in(0, 1]},\\
\beta & {\mbox{\scriptsize\normalsize if  } y=0, x\in(0, 1]},
\end{array}
\right.
\end{equation}
$\psi(x)=\left\{\begin{array}{ll}
0 & {\mbox{\scriptsize\normalsize if } x=0,}\\
1 & {\mbox{\scriptsize\normalsize if } x \in(0, 1]}
\end{array}
\right.$
and $\varphi:[0, 1] \rightarrow[0, b]$ is an increasing bijection.

\item [(\romannumeral 3)] There exist two increasing functions $h, g : [0, 1] \rightarrow[0, 1]$ fulfilling that $h(1)=g(1)=1$ such that for all $x, y\in [0,1]$
\begin{equation}\label{eq3.12}
A(x,y)=\left\{\begin{array}{ll}
1 & {\mbox{\scriptsize\normalsize if  } (x,y)=(1,1)},\\
0 & {\mbox{\scriptsize\normalsize if  } x,y \in[0, 1)},\\
g(y) & {\mbox{\scriptsize\normalsize if  } x=1, y\in[0, 1)},\\
h(x) & {\mbox{\scriptsize\normalsize if  } y=1, x\in[0, 1)},
\end{array}
\right.
\end{equation}
$\psi(x)=\left\{\begin{array}{ll}
0 & {\mbox{\scriptsize\normalsize if } x\in[0, 1),}\\
1 & {\mbox{\scriptsize\normalsize if } x=1}
\end{array}
\right.$
and $\varphi:[0, 1] \rightarrow[0, b]$ is an increasing bijection.
\end{enumerate}
\end{theorem}
\begin{proof}
Suppose that the aggregation function $A$ is $(\varphi, \psi)$ \emph{-quasi-homogeneous}. Then, by Lemmas \ref{lem3.2}, \ref{lemma3.1} and Remark \ref{Remark3.1}, we have $\psi(x)=x, \varphi(x)=\delta_{A}^{-1}(x)$ or
$$\psi(x)=\left\{\begin{array}{ll}
0 & {\mbox{\scriptsize\normalsize if } x=0,}\\
1 & {\mbox{\scriptsize\normalsize if } x \in(0, 1],}
\end{array}
\right.
\mbox{ or }
\psi(x)=\left\{\begin{array}{ll}
0 & {\mbox{\scriptsize\normalsize if  } x \in[0, 1),}\\
1 & {\mbox{\scriptsize\normalsize if  } x=1.}
\end{array}
\right.$$
Next, we investigate the representation of quasi-homogeneous aggregation function $A$ by distinguishing three cases.\\

Case 1. If $\psi(x)=x$ and $\varphi(x)=\delta_{A}^{-1}(x)$ for any $x\in [0, 1]$, then from Eq.\eqref{eq2.1},
$$A(x,y)=\left\{\begin{array}{ll}
0 & {\mbox{\scriptsize\normalsize if  } (x,y)=(0,0)},\\
\delta_{A}(y\delta_{A}^{-1}(A(\frac{x}{y},1))) & {\mbox{\scriptsize\normalsize if  } x\leq y \mbox{ and }y\neq0},\\
\delta_{A}(x\delta_{A}^{-1}(A(1,\frac{y}{x}))) & {\mbox{\scriptsize\normalsize if  } y\leq x \mbox{ and }x\neq0}.
\end{array}
\right.$$
Consider the functions $g,h:[0, 1] \rightarrow[0, 1]$ defined by $g(x)=A(1,x)$ and $h(x)=A(x,1)$, respectively. Then $g(1)=h(1)=1$, $g$ and $h$ are increasing. Moreover, $$A(x,y)=\left\{\begin{array}{ll}
0 & {\mbox{\scriptsize\normalsize if  } (x,y)=(0,0)},\\
\delta_{A}(y\delta_{A}^{-1}(h(\frac{x}{y}))) & {\mbox{\scriptsize\normalsize if  } x\leq y \mbox{ and }y\neq0},\\
\delta_{A}(x\delta_{A}^{-1}(g(\frac{y}{x}))) & {\mbox{\scriptsize\normalsize if  } y\leq x \mbox{ and }x\neq0}.
\end{array}
\right.$$
Now, we show that $\frac {\delta_{A}^{-1}(h(x))}{x}$ is decreasing for all $x\in(0, 1]$. Suppose $\frac {\delta_{A}^{-1}(h(x))}{x}$ is not decreasing for all $x\in(0, 1]$. Then there exist  $0<y_{1}<z_{1}\leq 1 $ such that $\frac {\delta_{A}^{-1}(h(y_{1}))}{y_{1}}<\frac {\delta_{A}^{-1}(h(z_{1}))}{z_{1}}$. Let $0<x\leq y<z\leq 1$ with $\frac {x}{y}=z_{1}$ and $\frac {x}{z}=y_{1}$. Then one has that $\frac {\delta_{A}^{-1}(h(\frac {x}{z}))}{\frac {x}{z}}<\frac {\delta_{A}^{-1}(h(\frac {x}{y}))}{\frac {x}{y}}$, which leads to $x\frac {\delta_{A}^{-1}(h(\frac {x}{z}))}{\frac {x}{z}}<x\frac {\delta_{A}^{-1}(h(\frac {x}{y}))}{\frac {x}{y}}$. Thus $z\delta_{A}^{-1}(h(\frac {x}{z}))<y\delta_{A}^{-1}(h(\frac {x}{y}))$.
So, for $0<x\leq y<z\leq 1$, we have that
$$z\delta_{A}^{-1}(A(\frac {x}{z},1))<y\delta_{A}^{-1}(A(\frac {x}{y},1)).$$
This follows that $\delta_{A}(z\delta_{A}^{-1}(A(\frac {x}{z},1)))<\delta_{A}(y\delta_{A}^{-1}(A(\frac {x}{y},1)))$, i.e., $A(x,z)<A(x,y)$, a contradiction. Therefore, $\frac {\delta_{A}^{-1}(h(x))}{x}$ is decreasing for all $x\in(0, 1]$. Similarly, we can prove that $\frac {\delta_{A}^{-1}(g(x))}{x}$ is decreasing for all $x\in(0, 1]$.

Case 2. If $$\psi(x)=\left\{\begin{array}{ll}
0 & {\mbox{\scriptsize\normalsize if } x=0,}\\
1 & {\mbox{\scriptsize\normalsize if } x \in(0, 1],}
\end{array}
\right.$$
then from Eq.\eqref{eq2.1}, for an arbitrary increasing bijection $\varphi:[0, 1] \rightarrow[0, b]$,
$$A(x, x) =\varphi^{-1}(\psi(x)\varphi(A(1, 1)))=A(1, 1)=1, x \in(0, 1].$$
So
\begin{equation*}
A(x,x)=\left\{\begin{array}{ll}
0 & {\mbox{\scriptsize\normalsize if  } x=0},\\
1 & {\mbox{\scriptsize\normalsize if  } x\in(0, 1]}
\end{array}
\right.
\end{equation*}
and
$$A(x,y)=1 \mbox{ for all } x,y \in(0, 1].$$
Moreover, $A(0,y)=A(y\cdot0,y\cdot1)=\varphi^{-1}(\psi(y)\varphi(A(0, 1)))=A(0, 1)$ for any $y\in(0, 1]$ and $A(x,0)=A(x\cdot1,x\cdot0)=\varphi^{-1}(\psi(x)\varphi(A(1, 0)))=A(1, 0)$ for any $x\in(0, 1]$. Let $\alpha=A(0, 1)$ $\beta=A(1, 0)$. Therefore,
\begin{equation*}
A(x,y)=\left\{\begin{array}{ll}
0 & {\mbox{\scriptsize\normalsize if  } (x,y)=(0,0)},\\
1 & {\mbox{\scriptsize\normalsize if  } x,y \in(0, 1]},\\
\alpha & {\mbox{\scriptsize\normalsize if  } x=0, y\in(0, 1]},\\
\beta & {\mbox{\scriptsize\normalsize if  } y=0, x\in(0, 1]}.
\end{array}
\right.
\end{equation*}

Case 3. If $$\psi(x)=\left\{\begin{array}{ll}
0 & {\mbox{\scriptsize\normalsize if  } x \in[0, 1),}\\
1 & {\mbox{\scriptsize\normalsize if  } x=1,}
\end{array}
\right.$$  Then from Eq.\eqref{eq2.1}, for an arbitrary increasing bijection $\varphi:[0, 1] \rightarrow[0, b]$,
$$A(x, x) =\varphi^{-1}(\psi(x)\varphi(A(1, 1)))=\varphi^{-1}(0)=0, x \in[0, 1).$$
So
\begin{equation*}
A(x,x)=\left\{\begin{array}{ll}
0 & {\mbox{\scriptsize\normalsize if  }x\in[0, 1)},\\
1 & {\mbox{\scriptsize\normalsize if  }  x=1}
\end{array}
\right.
\end{equation*}
and
$$A(x,y)=0 \mbox{ for all } x,y \in[0, 1).$$
Consider the functions $g,h:[0, 1] \rightarrow[0, 1]$ defined by $g(x)=A(1,x)$ and $h(x)=A(x,1)$, respectively. Therefore,
\begin{equation*}
A(x,y)=\left\{\begin{array}{ll}
1 & {\mbox{\scriptsize\normalsize if  } (x,y)=(1,1)},\\
0 & {\mbox{\scriptsize\normalsize if  } x,y \in[0, 1)},\\
g(y) & {\mbox{\scriptsize\normalsize if  } x=1, y\in[0, 1)},\\
h(x) & {\mbox{\scriptsize\normalsize if  } y=1, x\in[0, 1)}.
\end{array}
\right.
\end{equation*}

Conversely, one can directly verify that the aggregation function $A$ is $(\varphi, \psi)$-quasi-homogeneous if one of $(\romannumeral 1), (\romannumeral 2), (\romannumeral 3) $ holds.
\end{proof}

Notice that Theorem \ref{th3.2} describes the $(\varphi, \psi)$-quasi-homogeneity of aggregation function $A$ in the case $0<b<\infty$ completely. Although the inverse function $\varphi^{-1}$ of function $\varphi$ in the definition of $(\varphi, \psi)$-quasi-homogeneity of aggregation function $A$ in \cite{Su2022} can be regarded as a function from $\varphi([0,1])$ to $[0,1]$, one can easily see that Theorem 6 of \cite{Su2022} just concerns the $(\varphi, \psi)$-quasi-homogeneity of aggregation function $A$ when the diagonal $\delta_{A}$ is continuous, and it does not characterize the $(\varphi, \psi)$-quasi-homogeneity of aggregation function $A$ when the diagonal $\delta_{A}$ is not continuous.
\subsection{The case $b=\infty$}
In what follows, it is tacitly assumed that $0\cdot \infty=\infty\cdot 0 =0$ and $\frac {1} {\infty}=0$.
\begin{theorem}\label{th2.4}
An aggregation function $A:[0, 1]^2\rightarrow[0, 1]$ is $(\varphi, \psi)$-quasi-homogeneous  if and only if one of the following statements holds:
\begin{enumerate}
\item [(\romannumeral 1)] There exist two constants $\alpha, \beta \in [0, 1]$ such that for all $x, y\in [0,1]$,
$$A(x,y)=\left\{\begin{array}{ll}
0 & {\mbox{\scriptsize\normalsize if  } (x,y)=(0,0)},\\
1 & {\mbox{\scriptsize\normalsize if  } x,y \in(0, 1]},\\
\alpha & {\mbox{\scriptsize\normalsize if  } x=0, y\in(0, 1]},\\
\beta & {\mbox{\scriptsize\normalsize if  } y=0, x\in(0, 1]},
\end{array}
\right.$$
$$\psi(x)=\left\{\begin{array}{ll}
0 & {\mbox{\scriptsize\normalsize if } x=0,}\\
1 & {\mbox{\scriptsize\normalsize if } x \in(0, 1]}
\end{array}
\right.$$
and $\varphi:[0, 1] \rightarrow[0, b]$ is an increasing bijection.
\item [(\romannumeral 2)] There exist two increasing functions $h, g : [0, 1] \rightarrow[0, 1]$ fulfilling that $h(1)=g(1)=1$ such that for all $x, y\in [0,1]$,
$$A(x,y)=\left\{\begin{array}{ll}
1 & {\mbox{\scriptsize\normalsize if  } (x,y)=(1,1)},\\
0 & {\mbox{\scriptsize\normalsize if  } x,y \in[0, 1)},\\
g(y) & {\mbox{\scriptsize\normalsize if  } x=1, y\in[0, 1)},\\
h(x) & {\mbox{\scriptsize\normalsize if  } y=1, x\in[0, 1)},
\end{array}
\right.$$
$$\psi(x)=\left\{\begin{array}{ll}
0 & {\mbox{\scriptsize\normalsize if } x\in[0, 1),}\\
1 & {\mbox{\scriptsize\normalsize if } x=1}
\end{array}
\right.$$
and $\varphi:[0, 1] \rightarrow[0, b]$ is an increasing bijection.
\end{enumerate}
\end{theorem}
\begin{proof}
Considering $x=1$ in Eq.(\ref{eq3.1}), we have
\begin{equation}\label{eq2.2.1}
\varphi(\delta_{A}(\lambda))=\psi(\lambda)\varphi(1)
\end{equation}
 for all $\lambda \in[0, 1]$.
 By Eqs.(\ref{eq3.1}) and (\ref{eq2.2.1}), we have
 \begin{equation*}
\varphi(1)\psi(\lambda x)=\varphi(1)\psi(\lambda)\psi(x)\mbox{ for all } \lambda, x\in[0, 1].
 \end{equation*}
 Then $\varphi(1)[\psi(\lambda x)-\psi(\lambda)\psi(x)]=0$ for all $\lambda, x\in[0, 1]$, which means that \begin{equation}\label{eq2.2.4}
\psi(\lambda x)=\psi(\lambda)\psi(x)\mbox{ for all } \lambda, x\in[0, 1].
 \end{equation}
Considering $\lambda=1$ in Eq.(\ref{eq2.2.1}), we have
 $\varphi(\delta_{A}(1))=\psi(1)\varphi(1)$, i.e., $\varphi(1)=\psi(1)\varphi(1)$. Then $\varphi(1)[1-\psi(1)]=0$. This follows
 \begin{equation}\label{eq2.2.2}
 \psi(1)=1.
 \end{equation}
Considering $\lambda=0$ in Eq.(\ref{eq2.2.1}), we have
 $\varphi(\delta_{A}(0))=\psi(0)\varphi(1)$, i.e., $\varphi(0)=\psi(0)\varphi(1)$. Then
\begin{equation}\label{eq2.2.3}
\psi(0)=0
 \end{equation}
since $\varphi(0)=0$ by Lemma \ref{lem3.1} (ii). Again, from Eq.(\ref{eq2.2.1}), we have
$\varphi(1)=\psi(\lambda)\varphi(1)\mbox{ for any } \lambda \in(0, 1)$, i.e., $\psi(\lambda)=1 \mbox{ for any } \lambda \in(0, 1)$ if $\delta_{A}(\lambda)=1$,
and
$\varphi(0)=\psi(\lambda)\varphi(1)\mbox{ for any } \lambda \in(0, 1)$, i.e., $\psi(\lambda)=0 \mbox{ for any } \lambda \in(0, 1)$ if $\delta_{A}(\lambda)=0$.
In the following, we just need to prove that $\delta_{A}(\lambda)=1$ or $\delta_{A}(\lambda)=0$ for any $\lambda\in
(0,1)$. In fact, if $\delta_{A}(\lambda)<1$ for a $\lambda \in(0, 1)$, then $\varphi(\delta_{A}(\lambda))<\infty$.
Thus by Eq. (\ref {eq2.2.1}), we have
\begin{equation*}\label{eq2.2.5}
\psi(\lambda)=0.
\end{equation*}
Therefore,
\begin{equation*}
\varphi(\delta_{A}(\lambda))=\psi(\lambda)\varphi(1)=0.
\end{equation*}
Then from Lemma \ref{lem3.1}, $\delta_{A}(\lambda)=0$ since $\varphi$ is a bijection. Consequently, we distinguish three cases as below.
\begin{enumerate}
\item [(\romannumeral 1)] If $\delta_{A}(\lambda)=1$ for all $\lambda\in(0, 1)$, then by Eq.(\ref{eq2.2.1}) we have
$\psi(\lambda)=1$ for all $\lambda\in(0, 1)$. From Eqs.(\ref{eq2.2.2}) and (\ref{eq2.2.3}),
we have that for all $x\in [0,1]$,
$$\psi(x)=\left\{\begin{array}{ll}
0 & {\mbox{\scriptsize\normalsize if } x=0,}\\
1 & {\mbox{\scriptsize\normalsize if } x \in(0, 1].}
\end{array}
\right.$$
Therefore, from Eq.\eqref{eq2.1}, for an arbitrary increasing bijection $\varphi:[0, 1] \rightarrow[0, b]$,
$$A(x, x) =\varphi^{-1}(\psi(x)\varphi(A(1, 1)))=A(1, 1)=1, x \in(0, 1].$$
So
\begin{equation*}
A(x,x)=\left\{\begin{array}{ll}
0 & {\mbox{\scriptsize\normalsize if  } x=0},\\
1 & {\mbox{\scriptsize\normalsize if  } x\in(0, 1]}
\end{array}
\right.
\end{equation*}
and
$$A(x,y)=1 \mbox{ for all } x,y \in(0, 1].$$
Moreover, $A(0,y)=A(y\cdot0,y\cdot1)=\varphi^{-1}(\psi(y)\varphi(A(0, 1)))=A(0, 1)$ for any $y\in(0, 1]$ and $A(x,0)=A(x\cdot1,x\cdot0)=\varphi^{-1}(\psi(x)\varphi(A(1, 0)))=A(1, 0)$ for any $x\in(0, 1]$. Let $\alpha=A(0, 1)$ $\beta=A(1, 0)$.
Therefore, we have that for all $x, y\in [0,1]$,
$$A(x,y)=\left\{\begin{array}{ll}
0 & {\mbox{\scriptsize\normalsize if  } (x,y)=(0,0)},\\
1 & {\mbox{\scriptsize\normalsize if  } x,y \in(0, 1]},\\
\alpha & {\mbox{\scriptsize\normalsize if  } x=0, y\in(0, 1]},\\
\beta & {\mbox{\scriptsize\normalsize if  } y=0, x\in(0, 1]}.
\end{array}
\right.$$

\item [(\romannumeral 2)] If $\delta_{A}(\lambda)=0$ for all $\lambda\in(0, 1)$, then by Eq.\ref{eq2.2.1}, $\psi(\lambda)=0$ for all $\lambda\in(0, 1)$.
Then by Eqs.(\ref{eq2.2.2}) and (\ref{eq2.2.3}),
$$\psi(x)=\left\{\begin{array}{ll}
0 & {\mbox{\scriptsize\normalsize if } x\in[0, 1),}\\
1 & {\mbox{\scriptsize\normalsize if } x=1.}
\end{array}
\right.$$
Then from Eq.\eqref{eq2.1}, for an arbitrary increasing bijection $\varphi:[0, 1] \rightarrow[0, b]$,
$$A(x, x) =\varphi^{-1}(\psi(x)\varphi(A(1, 1)))=\varphi^{-1}(0)=0, x \in[0, 1).$$
So
\begin{equation*}
A(x,x)=\left\{\begin{array}{ll}
0 & {\mbox{\scriptsize\normalsize if  }x\in[0, 1)},\\
1 & {\mbox{\scriptsize\normalsize if  }  x=1}
\end{array}
\right.
\end{equation*}
and
$$A(x,y)=0 \mbox{ for all } x,y \in[0, 1).$$
Consider the functions $g,h:[0, 1] \rightarrow[0, 1]$ defined by $g(x)=A(1,x)$ and $h(x)=A(x,1)$, respectively. Therefore,
\begin{equation*}
A(x,y)=\left\{\begin{array}{ll}
1 & {\mbox{\scriptsize\normalsize if  } (x,y)=(1,1)},\\
0 & {\mbox{\scriptsize\normalsize if  } x,y \in[0, 1)},\\
g(y) & {\mbox{\scriptsize\normalsize if  } x=1, y\in[0, 1)},\\
h(x) & {\mbox{\scriptsize\normalsize if  } y=1, x\in[0, 1)}.
\end{array}
\right.
\end{equation*}
\item [(\romannumeral 3)] If there exists an $e\in(0, 1)$ such that
$$\delta_{A}(\lambda)=0\mbox{ for all } \lambda\in(0, e)$$
and
$$\delta_{A}(\lambda)=1\mbox{ for all } \lambda\in(e, 1),$$
then by Eq.(\ref {eq2.2.1}), we have
\begin{equation}\label{eq2.002.3}\psi(\lambda)=0\mbox{ for all } \lambda\in(0, e)\end{equation}
and
$$\psi(\lambda)=1\mbox{ for all } \lambda\in(e, 1).$$
Therefore, by Eq.(\ref {eq2.2.4}), there exists a $\lambda\in(e, 1)$ such that
\begin{equation}\label{eq2.0020.3}\psi(\lambda^n)=(\psi(\lambda))^n=1\mbox{ for all } n\in N.\end{equation}
On the other hand, there exists an $m\in N$ such that $0<\lambda^m<e$ since $\lim\limits_{n\rightarrow\infty}\lambda^{n}=0$. So it follows from Eq.\eqref{eq2.002.3} that $\psi(\lambda^{m})=0$, contrary to Eq.\eqref{eq2.0020.3}.
\end{enumerate}
\end{proof}

Notice that applying Theorems \ref{th3.2} and \ref{th2.4}, from the point of view of the functional representation, the quasi-homogeneous aggregation functions are exactly classified into three classes as Eqs.\eqref{eq3.10}, \eqref{eq3.11} and \eqref{eq3.12}, respectively. One can check that the diagonal function of the first classified aggregation function is continuous and the others are non-continuous. Moreover, the constructions of the last two classes are very clearly. Therefore, our next attention is how to construct the first classified quasi-homogeneous aggregation functions.

\section{Constructions of quasi-homogeneous aggregation functions}
In this section, we first introduce a triple generator of quasi-homogeneous aggregation function, then investigate how to construct quasi-homogeneous aggregation functions by using their triple generators.

\begin{definition}\label{def5.1}
\emph{Let $f:[0, 1] \rightarrow[0, 1]$ be an increasing bijection and $g,h:[0, 1] \rightarrow[0, 1]$ be two increasing functions. If a function $A:[0, 1]^2\rightarrow[0, 1]$ given by
\begin{equation}\label{eq5.1}
A(x,y)= \left\{\begin{array}{ll}
0 & {\mbox{ if }\  x=y=0,}\\
f(yf^{-1}(h(\frac{x}{y}))) & {\mbox{ if }\  x\leq y \mbox{ and }y\neq0,}\\
f(xf^{-1}(g(\frac{y}{x}))) & {\mbox{ if }\  y\leq x \mbox{ and }x\neq0.}
\end{array}
\right.
\end{equation}
is a quasi-homogeneous aggregation function, then $(f, g, h)$ is called a} triple generator \emph{of quasi-homogeneous aggregation function $A$, and $A$ is said to be a quasi-homogeneous aggregation function generated by the triple $(f, g, h)$.}
\end{definition}

\begin{proposition}\label{prop5.1}
 Let $f:[0, 1] \rightarrow[0, 1]$ be an increasing bijection on $[0,1]$ and $g, h:[0, 1] \rightarrow[0, 1]$ be two increasing functions fulfilling the following conditions:
 \begin{enumerate}
\item [(1)] $h(1)=1$, $g(1)=1$;
\item [(2)] $\frac {f^{-1}(h(x))}{x}$ and $\frac {f^{-1}(g(x))}{x}$ are two decreasing functions on $(0, 1]$.
\end{enumerate}
Then the function $A:[0, 1]^2\rightarrow[0, 1]$ defined by Eq.(\ref{eq5.1}) is a quasi-homogeneous aggregation function, i.e., $(f, g, h)$ is the triple generator of the quasi-homogeneous aggregation function $A$. Moreover, the diagonal $\delta_{A}$ is an increasing bijection.
\end{proposition}
\begin{proof}
 Obviously, $A(0,0)=0, A(1,1)=1$.
Next, we prove the monotonicity of function $A$. First, we prove that $A(x,y)\leq A(x,z)$ whenever $y\leq z$.
Obviously, $A(0, y)\leq A(0, z)$ whenever $y\leq z$. Considering $x\in(0, 1]$ and $y\leq z$.
\begin{enumerate}
\item [$\diamond$] If $y\leq z\leq x$, then
\begin{equation*}
\begin{split}
A(x,y)&=f(xf^{-1}(A(1,\frac{y}{x})))\\
&=f(xf^{-1}(g(\frac{y}{x})))\\
&\leq f(xf^{-1}(g(\frac{z}{x})))  \mbox{ by the monotonicity  of } f\mbox{ and }g\\
&=A(x,z).
\end{split}
\end{equation*}
\item [$\diamond$]If $y\leq x\leq z$, then $\frac {f^{-1}(h(x))}{x}\geq f^{-1}(h(1))/1$, i.e., $f^{-1}(h(x))\geq x$ since $\frac {f^{-1}(h(x))}{x}$ is decreasing on $(0, 1]$. Thus
\begin{equation*}
\begin{split}
A(x,y)&=f(xf^{-1}(g(\frac{y}{x})))\\
&\leq f(x)\\
&=f(z\cdot\frac{x}{z})\\
&\leq f(zf^{-1}(h(\frac{x}{z})))\\
&=A(x,z).
\end{split}
\end{equation*}
\item [$\diamond$] If $x\leq y\leq z$, then
\begin{equation*}
\begin{split}
A(x,y)&=f(yf^{-1}(h(\frac{x}{y})))\\
&\leq f(zf^{-1}(h(\frac{x}{z}))) \mbox{  by the monotonicity of } \frac {f^{-1}(h(x))}{x}\\
&=A(x,z).
\end{split}
\end{equation*}
\end{enumerate}

Bring all observations together, we can see that $A(x,y)\leq A(x,z)$ whenever $y\leq z$. Similarly, we can prove that $A(y,x)\leq A(z,x)$ whenever $y\leq z$.

Secondly, we prove the quasi-homogeneity of $A$. Considering the functions $\psi(x)=x, \varphi(x)=f^{-1}(x)$ for all $x\in[0,1]$. Then by Eq.(\ref{eq5.1}) for all $x, y, \lambda\in[0,1]$,
\begin{equation*}
\begin{split}
A(\lambda x, \lambda y)&=\left\{\begin{array}{ll}
0 & {\mbox{ if }\  x=y=0  \mbox{ or}\  \lambda=0,}\\
f(\lambda yf^{-1}(h(\frac{\lambda x}{\lambda y}))) & {\mbox{ if }\  \lambda x\leq \lambda y \ \mbox{and}\  \lambda y\neq0,}\\
f(\lambda xf^{-1}(g(\frac{\lambda y}{\lambda x}))) & {\mbox{ if }\  \lambda y\leq \lambda x \ \mbox{and}\ \lambda x\neq0.}
\end{array}
\right. \\
&=\left\{\begin{array}{ll}
0 & {\mbox{ if }\  x=y=0  \mbox{ or}\  \lambda=0,}\\
f(\lambda yf^{-1}(h(\frac{x}{y}))) & {\mbox{ if }\ x\leq y \ \mbox{and}\ y\neq0,}\\
f(\lambda xf^{-1}(g(\frac{y}{x}))) & {\mbox{ if }\ y\leq x \ \mbox{and}\ x\neq0.}
\end{array}
\right.\\
&=\left\{\begin{array}{ll}
0 & {\mbox{ if }\  x=y=0  \mbox{ or}\  \lambda=0,}\\
f(\lambda f^{-1}(f(yf^{-1}(h(\frac{x}{y}))))) & {\mbox{ if }\ x\leq y \ \mbox{and}\ y\neq0,}\\
f(\lambda f^{-1}(f(xf^{-1}(g(\frac{y}{x}))))) & {\mbox{ if }\ y\leq x \ \mbox{and}\ x\neq0.}
\end{array}
\right.\\
&=f(\lambda f^{-1}( A( x, y)))\\
&=\varphi^{-1}(\psi(\lambda )\varphi( A( x, y))).
\end{split}
\end{equation*}
Therefore, $A$ is $(\varphi, \psi)$-quasi-homogeneous with $\psi(x)=x, \varphi(x)=f^{-1}(x)$ for all $x\in[0,1]$. It is easy to see that $\delta_{A}(x,x)=f(x)$ for all $x\in[0,1]$. Thus $\delta_{A}$ is an increasing bijection.
\end{proof}

\begin{example}\label{example5.1}
\emph{Consider the functions $g, h:[0, 1]\rightarrow[0, 1]$ and $f:[0, 1]\rightarrow[0, 1]$ defined by $g(x)=x, h(x)=\frac {2x}{1+x}$ and $f(x)=x$, respectively. It is easy to verify that $g, h$ and $f$ satisfy the conditions of Proposition \ref{prop5.1}. Then by Proposition \ref{prop5.1}, the function $A:[0, 1]^2\rightarrow[0, 1]$ with
$$A(x,y)=\left\{\begin{array}{ll}
\frac {2xy}{x+y} & {\mbox{ if }\ x\leq y \ \mbox{and}\ y\neq0,}\\
y & {\mbox{ if }\ y\leq x}
\end{array}
\right.
$$
for all $x,y\in[0, 1]$ is a quasi-homogeneous aggregation function whose triple generator is $(f,g,h)$.}
\end{example}

From Theorem \ref{th3.2}, one can immediately has the following proposition.
\begin{proposition}\label{prop5.2}
Let $A:[0, 1]^2 \rightarrow[0, 1]$ be a quasi-homogeneous aggregation function. Then $A$ has a triple generator $(f,g, h)$ with
 \begin{enumerate}
\item [(1)] $f(x)=\delta_{A}(x)$ for all $x\in[0, 1]$;
\item [(2)] $g(x)=A(1, x), h(x)=A(x, 1)$ for all $x\in[0, 1]$.
\end{enumerate}
\end{proposition}

\section{Conclusions}

 The contributions of this article include two aspects, one is that it characterized all quasi-homogeneous aggregation functions, which classifies the quasi-homogeneous aggregation functions into three classes. The other is that it constructed quasi-homogeneous aggregation functions by using their triple generators. All the results show us that generally we should mainly consider the quasi-homogeneous aggregation functions whose diagonal functions are continuous when we explore the quasi-homogeneity of aggregation functions. Meanwhile, it should be stressed that a slightly modifying Definition \ref{def2.2} does not influence the currently existing results in literature.

\end{document}